\documentclass[12pt,psfig]{article}
\usepackage{graphicx,epsfig}
\newcommand{\si}{\sigma}
\newcommand{\ld}{\lambda}
\newcommand{\ba}{\begin{array}}
\newcommand{\ea}{\end{array}}
\newcommand{\ban}{\begin{eqnarray*}}
\newcommand{\ean}{\end{eqnarray*}}

\newcommand{\la}{\langle}
\newcommand{\ra}{\rangle}

\begin{document}

\baselineskip=17pt

\begin{center}

\vspace{-0.6 in} {\large \bf  Stability of essential spectra\\
 of self-adjoint subspaces under compact perturbations}

 \vspace{0.3in}

Yuming Shi \\

Department of Mathematics, Shandong University\\

Jinan, Shandong 250100, P. R. China\\

Email address: ymshi@sdu.edu.cn

\end{center}

{\bf Abstract.} This paper studies stability of essential spectra
of self-adjoint subspaces (i.e., self-adjoint linear relations)
under finite rank and compact perturbations in Hilbert spaces.
Relationships between compact perturbation of closed subspaces
and relatively compact perturbation of their operator parts
are first established. This gives a characterization of
compact perturbation in terms of difference
between the operator parts of perturbed and unperturbed subspaces.
It is shown that a self-adjoint subspace is still self-adjoint
under either relatively bounded perturbation with relative bound less than one
or relatively compact perturbation or compact perturbation
with a certain additional condition. By using these results,
invariance of essential spectra of self-adjoint subspaces
is proved under relatively compact and compact perturbations, separately.
As a special case, finite rank perturbation is discussed.
The results obtained in this paper generalize the corresponding results for self-adjoint operators
to self-adjoint subspaces.

\medskip

\noindent{\bf 2010 AMS Classification}: 47A06, 47A55, 47A10, 47B25.

\bigskip

\noindent{\bf 1. Introduction}\medskip

Perturbation problems are one of the main topics in both pure
and applied mathematics.  The perturbation theory of operators
(i.e., single-valued operators) has been extensively studied
and many elegant results have been obtained (cf., [7, 9, 16]).
In particular, stability of spectra of
self-adjoint operators under perturbation has got lots of attention.
We shall recall several most well-known results about it.
If a perturbation term is a symmetric (i.e., densely-defined
and Hermitian) and relatively compact operator to a self-adjoint operator,
then the essential spectrum of the self-adjoint operator
 is invariant  (see [11, Theorem 8.15] or [16, Theorem 9.9]).
However, it was shown that
its absolutely continuous spectrum may disappear
under this perturbation even though the
perturbation term is very small by H. Weyl
and later generalized by von Neumann [7, Chapter 10, Theorem 2.1].
But if the perturbation term  is finite rank or
more generally belongs to the trace class,
then its absolutely continuous spectrum is invariant
 [7, Chapter 10, Theorems 4.3 and 4.4].
 These results have been extensively applied to study of
 stability of spectra of symmetric linear differential operators
 and bounded Jacobi operators (i.e, second-order  bounded and
 symmetric linear difference operators)
 including Schr\"odinger operators that have a strong physical
 background.

Recently, it was found that minimal and maximal operators generated by symmetric
linear difference expressions are multi-valued or non-densely defined
in general even though the corresponding definiteness condition is satisfied (cf.,
[10, 14]), and similar are those generated
by symmetric linear differential expressions that do not
satisfy the definiteness condition [8].
So the classical perturbation theory of operators
are not available in this case.
Partially due to the above reason, the study of non-densely defined
or multi-valued operators has attracted a great deal of interests
in near half a century.

In 2009, Azizov with his coauthors introduced concepts of
compact and finite rank perturbations of closed subspaces
in $X\times Y$ in terms of difference between
orthogonal projections of $X\times Y$
to the subspaces, where $X$ and $Y$ are Hilbert spaces [2].
They proved that a closed subspace is a finite rank or compact
perturbation of another closed subspace if and only if
the difference between their resolvents is a finite rank
or compact operator in the intersection of their
resolvent sets in the case that this intersection is not empty
and $X=Y$ [2, Corollaries 3.4 and 4.5]. Further, they studied
stability properties of spectral points of positive and
negative type and type $\pi$ in the non-self-adjoint
case under several kinds of perturbations in the Krein spaces [3].

Minimal operators or subspaces, generated by symmetric
linear differential and difference expressions, are
closed Hermitian operators or subspaces, and their self-adjoint
extensions are self-adjoint operators or subspaces.
Their resolvents can be expressed by
corresponding Green functions. In the case that the
 differential and difference expressions are singular,
their resolvents are complicated in general, and
much more complicated when their orders (or dimensions)
are higher because several boundary conditions
or coupled boundary conditions are involved.
Moreover, the Green functions are often expressed by
solutions of the systems rather than by coefficients of
the systems. Therefore,  in some cases it is more convenient to give
out a characterization of perturbation
in terms of the operators or subspaces themselves rather than
their resolvents.
For the operator case, concepts of relatively compact,
and finite rank and trace class perturbations were given
in terms of the difference between perturbed and
unperturbed operators (see Definition 2.3 and Lemma 2.4 for
relatively compact perturbation of operators in Section 2,
and we refer to [7, 11, 12, 16] for more detailed discussions).
We shall try to give a similar characterization of compact
and finite rank perturbations for closed subspaces to that for operators
in the present paper.

In this paper, we focus on the stability of essential spectra of self-adjoint
subspaces (i.e., self-adjoint linear relations)
under perturbations in Hilbert spaces.
Note that the spectrum and various spectra of a self-adjoint subspace,
including point, discrete, essential, continuous,
singular continuous, absolutely continuous and singular spectra,
can be only determined by its operator part [13, Theorems 2.1, 2.2,
3.4 and 4.1]. So it is natural to take into account perturbations
of operator parts of the unperturbed and perturbed subspaces.
However, we find that a summation of closed subspaces, $T=S+A$,
can not imply a similar relation of their operator
parts, $T_s=S_s+A_s$, in general
(see Section 3 for a detailed discussion).
In addition, in dealing with the summation of subspaces, one shall encounter
another problem that does not happen for operators.
If $T$, $S$ and $A$ are operators with $D(T)=D(S)\subset D(A)$ and satisfy
$T=S+A$, then $T$ can be interchanged with $S$ or $A$ as
$S=T-A$ or $A|_{D(S)}=T-S$. However, if they are multi-valued,
then this interchanging may not hold in general.
This is resulted in by their multi-valued parts.
See Example 3.1, in which $T=S+A$ holds, but $A|_{D(S)}\ne T-S$.
These problems make the study of subspaces complicated.

The rest of this paper is organized as follows.
In Section 2, some notations, basic concepts
and fundamental results about subspaces are introduced.
Relatively bounded and compact perturbations and finite rank perturbation of
operators and closed subspaces with some properties are recalled.
In Section 3,  relationships among the operator
parts of an  unperturbed subspace, its perturbation and its
perturbation term are established. Due to these relationships,
relationships of compact and finite rank perturbations of closed subspaces
with relatively compact and finite rank perturbations of the operator parts
of the unperturbed subspace and its perturbation term are given,
separately, in which characterizations of compact and finite rank
 perturbations are provided with relatively compact and finite rank
perturbations of their operator parts, respectively. In Section 4,
it is shown that a self-adjoint subspace is still
self-adjoint under either relatively bounded perturbation with
relative bound less than one or relatively compact perturbation
or compact perturbation with an additional condition.
Finally, it is proved that essential spectrum of a self-adjoint subspace
is stable under either compact perturbation or relatively compact
perturbation in Section 5.

\bigskip

\noindent{\bf 2. Preliminaries}\medskip

In this section, we shall first list some notations and basic concepts,
including spectrum,  discrete and essential spectra of subspaces,
and reducing subspace. Then we recall some fundamental results about subspaces.
Next, we introduce concepts of compact and finite rank perturbations for
operators and closed subspaces, and list some related results.

By ${\mathbf R}$ and ${\mathbf C}$ denote the sets of the real and
complex numbers, respectively, throughout this paper.

Let $X$ be a complex Hilbert space
with inner product $\langle \cdot, \cdot\rangle$,
and $T$ a linear subspace (briefly, subspace)
in the product space $X^2:=X\times
X$ with the following induced inner product, still denoted
by $\langle \cdot, \cdot\rangle$ without any confusion:
$$\langle(x,f),(y,g)\rangle=\langle x, y\rangle
+\langle f,g\rangle,\;\;(x,f),\;(y,g)\in X^2.$$
By $D(T)$ and $R(T)$ denote the domain and range of $T$,
respectively. Its adjoint subspace is defined by
$$T^*:=\{(y,g)\in X^2: \;\langle g,
x\rangle=\langle y,f\rangle\; {\rm  for\; all}\; (x, f)\in T\}.$$
$T$ is said to be an Hermitian subspace in $X^2$ if
$T\subset T^*$, and said to be a self-adjoint subspace in $X^2$ if
$T= T^*$.

Further, denote
$$T(x):=\{ f\in X:\; (x, f)\in T\}, \;\;
T^{-1}:=\{ (f,x):\; (x, f)\in T\}.$$
It is evident that $T(0) = \{0\}$ if and only if
$T$ uniquely determines a linear operator
from $D(T)$ into $X$ whose graph is $T$. For convenience,
a linear operator in $X$ will always be
identified with a subspace in $X^2$ via its graph.

Let $T$ and $S$ be two subspaces in $X^2$ and $\alpha\in {\bf C}$.
 Define
$$\begin{array}{ccc}\alpha T :
= \{(x, \alpha f) : (x, f) \in T\},\\[0.5ex]
T + S := \{(x, f + g) : (x, f) \in T,\; (x, g)\in S\}.
\end{array}$$
It is evident that if $T$ is closed, then $T- \lambda I$
is closed for any $\lambda\in {\bf C}$ and
$$(T-\lambda I)^*=T^*-\bar{\lambda}I.                                    \eqno(2.1) $$
On the other hand, if $T \cap S$ = \{(0, 0)\}, then denote
$$T\dot{+} S:=\{(x+y,f+g): \;(x,f)\in T,\; (y,g)\in S\}.$$
Further, if $T$ and $S$ are orthogonal; that is,
$\langle(x, f),(y, g)\rangle = 0 $ for all
$(x, f) \in T$ and $(y, g) \in S$, then we set
$$T\oplus S:=T\dot{+} S.$$

The following concepts were introduced in [5, 6, 13].\medskip

 \noindent{\bf Definition 2.1.} Let $T$ be a subspace in $X^2$.
\begin{itemize}
\item[{\rm (1)}] The set $\rho(T):=\{\ld\in {\mathbf C}:\,
(\ld I-T)^{-1}\;{\rm is\;  a\; bounded\;
 linear\; operator\; defined\; on}\;X\}$
 is called the resolvent set of $T$.

\item[{\rm (2)}]The set $\si(T):={\mathbf
C}\setminus \rho(T)$ is called the spectrum of $T$.

\item[{\rm (3)}] The essential spectrum $\sigma_e(T)$ of $T$
is the set of those points of $\si(T)$ that are either accumulation points of
$\sigma (T)$ or isolated eigenvalues of infinite multiplicity.

\item[{\rm (4)}] The set $\sigma_d(T):=\sigma (T)\setminus
\sigma_e(T) $ is called the discrete spectrum of $T$.
\end{itemize}

Arens [1] introduced the following
decomposition for a closed subspace $T$ in $X^2$:
$$T=T_s\oplus T_\infty,                                        \eqno (2.2)$$
where
$$T_\infty :=\{(0,g)\in X^2: (0,g)\in T\},\;\;
T_s :=T\ominus T_{\infty}.$$
Then $T_s$ is an operator. So $T_s$ and $T_\infty$ are called the operator
and pure multi-valued parts of $T$, respectively.
This decomposition will play an important role in our study.
Now, we shall recall their fundamental properties. The following
come from [1]:
$$D(T_s)=D(T),\;R(T_s)\subset T(0)^\perp,\;
R(T_{\infty})=T(0),\;\;T_\infty =\{0\}\times T(0),                               \eqno (2.3)$$
and $D(T_s)$ is dense in $T^*(0)^\perp$.

Throughout the present paper, the resolvent set and spectrum
of $T_s$ and $T_\infty$ mean those of
$T_s$ and $T_\infty$ restricted to $(T(0)^\perp)^2$ and $T(0)^2$,
respectively.
\medskip

\noindent{\bf Lemma 2.1 {\rm [13, Proposition 2.1 and Theorems 2.1, 2.2
and 3.4]}.}
Let $T$ be a closed Hermitian subspace in $X^2$. Then
$$T_s=T\cap (T(0)^\perp)^2,\;\;T_\infty=T\cap T(0)^2,                    \eqno (2.4)$$
$T_s$ is a closed Hermitian operator in $T(0)^\perp$,
$T_\infty$ is a closed Hermitian subspace in $T(0)^2$, and
$$\begin{array}{cc}
\rho(T)=\rho(T_s),\; \si(T)=\si(T_s),\; \si(T_\infty)=\emptyset,\\[0.8ex]
\si_p(T)=\si_p(T_s),\;\si_e(T)=\si_e(T_s),\;\si_d(T)=\si_d(T_s).
\end{array}                                                                            $$

In [4], Dijksma and Snoo introduced the concept of reducing subspace for a subspace.
Let $T$ be a subspace in $X^2$, $X_1$ a closed subspace in $X$ and  $P: X\to X_1$ an orthogonal projection. Denote
$$P^{(2)}T:=\{(Px,Pf):\,(x,f)\in T\}.$$
It is clear that $T\cap X_1^2\subset P^{(2)}T$.
If $P^{(2)}T\subset T$, then $X_1$ is called a reducing subspace of $T$.
We also say that $X_1$ reduces $T$ or $T$ is reduced by $X_1$.
In this case one has
$$T\cap X_1^2=P^{(2)}T.                                                  $$

\noindent{\bf Lemma 2.2 {\rm [4, Page 26]}.}
If $T$ is a self-adjoint subspace in $X^2$, then $T$ is reduced by $T(0)$, and $T_\infty$ and $T_s$
are self-adjoint subspaces in $T(0)^2$ and $(T(0)^\perp)^2$,
respectively.
\medskip

\noindent{\bf Lemma 2.3 {\rm [13, Theorem 2.5]}}.
 Let $T$ be an Hermitian subspace in $X^2$.
 Then $T$ is self-adjoint in $X^2$ if and only if
 $R(T-\ld I)=R(T-{\bar \ld I})=X$ for some $\ld \in {\mathbf C}$.

\medskip

Now, we present the following simple and useful result:\medskip

\noindent{\bf Proposition 2.1}.
 Let $T$ be a closed Hermitian subspace in $X^2$.
 Then $T$ is a self-adjoint subspace in $X^2$ if and only if $T_s$ is a self-adjoint
 operator in $T(0)^\perp$.\medskip

\noindent{\bf Proof.} The necessity directly follows from (ii) of Lemma 2.2.
Now, we consider the sufficiency. Since  $T_s$ is a self-adjoint
 operator in $T(0)^\perp$, one has that  $R(T_s\pm i I)=T(0)^\perp$
 by Theorem 5.21 of [16] or by Lemma 2.3.
So it follows from (2.2)-(2.4) that
$$R(T\pm i I)=R(T_s\pm i I)\oplus R(T_\infty)
 =T(0)^\perp\oplus T(0)=X,                                              $$
 which implies that $T$ is self-adjoint in $X^2$ by Lemma 2.3.
 The proof is complete.

\medskip

Next, we recall the concepts of relatively bounded
and compact operators and a related result, which will be used
in the sequels.\medskip

\noindent{\bf Definition 2.3 {\rm [16, Pages 93 and 275]}}.
Let $X$ be a Hilbert space, and $S$ and $T$ operators
in $X$. By $\|\cdot\|_S$
denote the graph norm of $S$, i.e.,
$\|x\|_S=\|x\|+\|Sx\|,\;\;x\in D(S).$
\begin{itemize}
\item [{\rm (1)}] $T$ is said to be $S$-bounded if $D(S)\subset D(T)$ and
there exists a constant $c\ge 0$ such that
$$\|Tx\|\le c\|x\|_S,\;\;x\in D(S).                                                 $$
\item [{\rm (2)}] If $T$ is $S$-bounded, then the infimum of all numbers $a\ge 0$
for which a constant $b\ge 0$ exists such that
$$\|Tx|\le a\|Sx\|+b \|x\|,\;\;\forall\;x\in D(S),      $$
is called the $S$-bound of $T$.
\item [{\rm (3)}] $T$ is said to be $S$-compact if it is compact as a mapping
from $(D(S),\|\cdot\|_S)$ into $X$.
\end{itemize}
\medskip

The following result is classical in the perturbation theory of
self-adjoint operators (cf., [11, Theorem 8.15] or
[16, Theorem 9.9]).\medskip

\noindent{\bf Lemma 2.4.} Let $T$ be a self-adjoint
operator in $X$, and $V$ a densely defined Hermitian and $T$-compact operator.
Then $T+V$ is a self-adjoint operator and $\si_e(T+V)=\si_e(T)$.
\medskip

To end this section, we recall concepts of finite rank
and compact perturbations for closed subspaces in the special case that
domains and ranges of the subspaces lie in a same Hilbert space.
We refer to [2] for more general definitions. \medskip

\noindent{\bf Definition 2.4.} Let $T$ and $S$ be closed subspaces
in $X^2$, and
$$P_T:\,X^2\to T,\;\;P_S:\,X^2\to S                                                    $$
 orthogonal projections.
\begin{itemize}
\item[{\rm (1)}] $T$ is said to be a compact perturbation of $S$
if $P_T-P_S$ is a compact operator in $X$.

\item[{\rm (2)}] $T$ is said to be a finite rank perturbation
of $S$ if $P_T-P_S$ is a finite rank operator in $X$.
\end{itemize}
\medskip

The following result gives characterizations
of compact and finite rank perturbations
of closed subspaces in terms of resolvents.\medskip

\noindent{\bf Lemma 2.5 {\rm [2, Corollaries 3.4 and 4.5]}.}
Let $T$ and $S$ be closed subspaces in $X^2$ and $\rho(T)\cap\rho(S) \ne \emptyset$.
Then
\begin{itemize}
\item[{\rm (i)}] $T$ is a compact perturbation of $S$
if and only if $(T-\ld I)^{-1}-(S-\ld I)^{-1}$ is a compact operator
in $X$ for some (and hence for all) $\ld \in \rho(T)\cap\rho(S)$;

\item[{\rm (ii)}] $T$ is a finite rank perturbation of $S$
if and only if $(T-\ld I)^{-1}-(S-\ld I)^{-1}$ is a finite
rank operator in $X$
for some (and hence for all) $\ld \in \rho(T)\cap\rho(S)$.
\end{itemize}

\bigskip

\noindent{\bf 3. Relationships among operator parts
and characterizations of compact and finite rank perturbations
}\medskip

In this section, we study relationships among operator parts
of an unperturbed subspace, its perturbation and its perturbation term.
Using them, we characterize compact and finite rank perturbations
of a closed subspace in terms of difference between their operator parts.

 we first study relationships
between operator parts of an unperturbed subspace, its perturbation and its perturbed term. Let $T$, $S$ and $A$ be closed subspaces in $X^2$
with $D(T)=D(S)\subset D(A)$, and
satisfy
$$T=S+A,                                                                     \eqno(3.1)$$
where $T$ can be regarded as a perturbation of $S$ by
the perturbation term $A$. It is natural
to ask whether their operator parts satisfy
$$T_s=S_s+A_s.                                                            \eqno(3.2)$$
This is a very interesting question itself.\medskip

The following simple fact will be repeatedly used in the sequent
discussion. If $S$ and $A$ are Hermitian subspaces, then
$T$, defined by (3.1), is an Hermitian subspace in $X^2$.
But if $S$ and $A$ are closed, $T$ may not be closed in general. \medskip

We first consider the following example:\medskip

\noindent{\bf Example 3.1.} Let $X=l^2[1,\infty)$, $e_1=\{e_{1n}\}_{n=1}^\infty\in X$
with $e_{11}=1$ and $e_{1n}=0$ for all $n\ge 2$, and $X_1={\rm span}\{e_1\}^\perp
=\{x=\{x_n\}_{n=1}^\infty\in X:\;x_1=0\}$. Define
$$S=\{(x,x+ce_1):\;x\in X_1, c\in {\mathbf C}\},                          \eqno(3.3)$$
and define $A$ by the following Jacobi operator:
$$(Ax)_n=a_{n}x_{n+1}+b_nx_n+a_{n-1}x_{n-1},\;n\ge 1,\; x\in X,              \eqno(3.4)$$
with $x_0=0$, where $a_n\ne 0$ for $n\ge 0$ and
$\{a_n,b_n\}$ is real and bounded.
Then $S$ is a closed Hermitian subspace with $D(S)=X_1$,
and $A$ is a bounded self-adjoint operator in $X$ with $D(A)=X$
(see [14] and [15] for more detailed discussions).

Let $T$ be defined  by (3.1). Then
$T$ is a closed Hermitian subspace in $X^2$, and
$$\begin{array}{cc}
T=\{(x,x+Ax+ce_1):\;x\in X_1,c\in {\mathbf C}\},
\\[0.6ex]
D(T)=D(S)=X_1\subset D(A),\;T(0)=S(0)={\rm span}\{e_1\}.
\end{array}$$
So $T(0)^\perp=S(0)^\perp=X_1.$ Further, by Lemma 2.1 one has
$$T_s=T\cap (T(0)^\perp)^2=T\cap X_1^2.$$
Thus, $(x,x+Ax+ce_1)\in T_s$ if and only if $x\in X_1$
and $x+Ax+ce_1\in X_1$. Since $x\in D(T)=X_1$, one has that
$(x,x+Ax+ce_1)\in T_s$ if and only if
 $Ax+ce_1\in X_1$, which is equivalent to $(Ax)_1+c=0$; that
is, $a_1x_2+c=0$ by (3.4), which yields
$c=-a_1x_2$. Hence,
$$T_sx=x+Ax-a_1x_2e_1,\;\;x\in X_1.                                                \eqno(3.5) $$
In addition, we can get by (3.3) and Lemma 2.1 that
$$S_sx=x,\;\;x\in X_1.                                                           \eqno(3.6) $$
It follows from (3.4)-(3.6) that
$$T_sx\ne S_sx+Ax,\;x\in X_1\;{\rm with}\;x_2\ne 0.                                   \eqno(3.7)$$
This means that (3.2) does not hold in general.
In particular, $T_s$ and $S_s+A_s$ have no inclusion relationships
in this example.
Furthermore, it follows from (3.4)-(3.6) that they satisfy
$$T_sx=S_sx+A_sx+Bx,\;\;x\in X_1,$$
where $Bx=-a_1\la e_2,x\ra e_1$ and $e_2=\{e_{2n}\}_{n=1}^\infty$
with $e_{22}=1$ and $e_{2n}=0$ for $n\ne 2$.
Obviously,  $B$ is a rank one operator.
\medskip

The following result gives a general relationship among $T_s$, $S_s$
and $A_s$.\medskip

\noindent{\bf Theorem 3.1.} Let $T$, $S$ and $A$ be closed
subspaces in $X^2$ with $D(T)=D(S)=:D\subset D(A)$, and satisfy (3.1).
And let
$$P:\;S(0)^\perp\to T(0)^\perp,\;\;Q:\;A(0)^\perp\to T(0)^\perp,            \eqno(3.8)$$
be orthogonal projections. Then
$$T_sx=PS_sx+QA_sx,\;\;x\in D.                                            \eqno(3.9)$$
Furthermore, if $S$ and $A$ are Hermitian subspaces in $X^2$,
then $PS_s$ and $QA_s$ are Hermitian operators in $D$, respectively.\medskip

\noindent{\bf Proof.} By (2.3) we have
$$\begin{array}{cc}
R(T_s)\subset T(0)^\perp, \;R(S_s)\subset S(0)^\perp,\;
R(A_s)\subset A(0)^\perp,\\[0.5ex]
D(T_s)=D(T)=D(S_s)=D\subset D(A_s)=D(A).
\end{array}                                                                     \eqno(3.10)    $$
It follows from (3.1) that
$$T(0)=S(0)+A(0),                                                          \eqno(3.11)$$
which is a sum of sets. This implies that
$S(0)\subset T(0)$ and $A(0)\subset T(0),$
and then
$$T(0)^\perp\subset S(0)^\perp, \; T(0)^\perp\subset A(0)^\perp.              \eqno(3.12)  $$
So the projections $P$ and $Q$ are well defined.

Fix any $x\in D$ and let $f=T_sx$. Then $f\in T(0)^\perp$
by (3.10) and $(x,f)\in T_s\subset T$ by (2.2).
So there exist $(x,g)\in S$ and $(x,h)\in A$ such that
$f=g+h$ by (3.1). Further, from (2.2) and (3.10),
there exist $g_s\in S(0)^\perp$, $g_\infty\in S(0)$,
$h_s\in A(0)^\perp$ and $h_\infty\in A(0)$ such that
$$
g=g_s+g_\infty, \;h=h_s+h_\infty,\;g_s=S_sx, \;h_s=A_sx.
                                                                        \eqno(3.13)$$
Then
$$\begin{array}{rrll}
&&f=g_s+h_s+g_\infty+h_\infty\\[0.5ex]
&=& Pg_s+Qh_s+(I-P)g_s+(I-Q)h_s+g_\infty+h_\infty.
\end{array}                                                                \eqno(3.14)   $$
Note that
$$(I-P)g_s\in S(0)^\perp\ominus T(0)^\perp\subset T(0),\;
(I-Q)h_s\in A(0)^\perp\ominus T(0)^\perp\subset T(0),                             $$
and $g_\infty, h_\infty\in T(0)$ by (3.11).
Hence, it follows from (3.14) that $(I-P)g_s+(I-Q_s)h_s+g_\infty+h_\infty=0$
and
$$f=Pg_s+Qh_s,$$
Therefore, (3.9) holds.

Further, suppose that $S$ and $A$ are Hermitian subspaces in $X^2$.
Then $T$ is an Hermitian subspace in $X^2$.
By (2.3) and Lemma 2.1, $T_s$, $S_s$ and
$A_s$ are closed Hermitian operators
in $T(0)^\perp$, $S(0)^{\perp}$ and $A(0)^\perp$, respectively, with
$$\begin{array}{cc}
T_s=T\cap(T(0)^\perp)^2,\;S_s=S\cap (S(0)^\perp)^2,\;
 A_s=A\cap (A(0)^\perp)^2,\\[0.8ex]
D(S_s)=D=D(T_s)\subset T(0)^\perp,\;\;
D(A_s)=D(A)\subset A(0)^\perp.
\end{array}                                                         \eqno(3.15)$$
Now, we show that $PS_s$ is  an Hermitian operator in $T(0)^\perp$.
For any given $x,y\in D$, let $S_sx=f_1+f_2$ and $S_sy=g_1+g_2$,
where $f_1,g_1\in T(0)^\perp$ and
$f_2,g_2\in S(0)^\perp\ominus T(0)^\perp\subset T(0)$.
Then $PS_sx=f_1$ and $PS_sy=g_1$. By noting
that $x,y\in D\subset T(0)^\perp$, it follows that
$\la f_2,y\ra=\la x,g_2\ra=0$. Hence, we have that
$$0=\la S_sx,y\ra-\la x,S_sy\ra =\la f_1,y\ra-\la x,g_1\ra=
\la PS_sx,y\ra-\la x,PS_sy\ra.$$
Therefore, $PS_s$ is an Hermitian operator.
With a similar argument, one can show that
$QA_s$ is an Hermitian operator in $D$.
This completes the proof.\medskip

It is evident that $T(0)=S(0)$ if $A$ is an operator in (3.1).
So the following result can be directly derived from Theorem 3.1.\medskip

\noindent{\bf Corollary 3.1.} Let $T$ and $S$ be closed
subspaces in $X^2$, $A$ is a closed operator in $X$
with $D(T)=D(S)=:D\subset D(A)$, and they satisfy (3.1). Then
$$T_sx=S_sx+QAx,\;\;x\in D,                                            $$
where $Q:\;X\to S(0)^\perp$ is an orthogonal projection.
\medskip

The unperturbed subspace $S$ is often self-adjoint
in $X^2$ in applications. Under this condition, we can get
a similar result to that in Corollary 3.1, which is
better than that in Theorem 3.1. \medskip

\noindent{\bf Theorem 3.2.} Let $T$ and $A$ be closed Hermitian
subspaces in $X^2$, $S$ a self-adjoint subspace in $X^2$
with $D(T)=D(S)=:D\subset D(A)$, and they satisfy (3.1).
Then
$$T_sx=S_sx+QA_sx,\;\;x\in D,                                            \eqno(3.16) $$
and
$$A(0)\subset S(0)=T(0),                                            \eqno(3.17) $$
where $Q:\, A(0)^\perp\to S(0)^\perp$ is an orthogonal projection.\medskip

\noindent{\bf Proof.} We first show that (3.17) holds.
By the assumptions, (3.9), (3.11), (3.12) and
(3.15) hold, and $S_s$ is a self-adjoint operator
in $S(0)^\perp$ by Lemma 2.2. So $D$ is dense in $S(0)^\perp$.
This yields that ${\bar D}=S(0)^\perp.$ In addition, by (3.15)
one has that $D\subset D(A)\subset A(0)^\perp$. Thus,
$S(0)^\perp\subset A(0)^\perp$, which, together
 with the fact that $A(0)$ and $S(0)$ are closed,
implies that $A(0)\subset S(0)$. Hence, (3.17) holds by (3.11).

Note that $P$, defined by (3.8), is an identity mapping from $S(0)^\perp$
onto itself in this case. Therefore, (3.16) follows.
The proof is complete.\medskip

The following result gives a sufficient condition such that
(3.2) holds.\medskip

\noindent{\bf Theorem 3.3.} Let $T$, $S$ and $A$ be closed Hermitian
subspaces in $X^2$ with $D(T)=D(S)=:D\subset D(A)$, and satisfy (3.1).
If $T(0)^\perp$ reduces $S$ and $A$, then
$$T_sx=S_sx+ A_sx,  \;\;x\in D,                                           \eqno(3.18) $$
$${\hat S_s}:=S\cap (T(0)^\perp)^2= S_s,
\;{\hat A_s}:=A\cap (T(0)^\perp)^2\subset A_s,                         \eqno(3.19)    $$
and ${\hat A_s}$ is a closed Hermitian operator in $T(0)^\perp$
 with $D\subset D({\hat A_s})$. \medskip

\noindent{\bf Proof.} By (2.3) and Lemma 2.1, $T_s$, $S_s$ and $A_s$
are closed Hermitian operators
in $T(0)^\perp$, $S(0)^{\perp}$ and $A(0)^\perp$, respectively,
and (3.11), (3.12) and (3.15) hold. From (3.12) and (3.15), we have
$${\hat S_s}\subset S_s,\; {\hat A_s}\subset A_s.                        \eqno(3.20)   $$
By noting that $T(0)^\perp$ is a closed subspace in $X$, ${\hat S_s}$ and ${\hat A_s}$
are closed Hermitian operators in $T(0)^\perp$.

We first show that the first relation in (3.19) holds.
By (3.20), it suffices to show that
$$S_s\subset {\hat S_s}.                                                  \eqno(3.21)$$
Fix any $(x,f)\in S_s$. Then $(x,f)\in (S(0)^\perp)^2$ and
 $x\in D\subset T(0)^\perp$ by (3.15).
There exist $f_1\in T(0)^\perp$ and
 $f_2\in S(0)^\perp\ominus T(0)^\perp$
such that
$$(x,f)=(x,f_1)+(0,f_2),                                            \eqno(3.22)  $$
where the first relation in (3.12) has been used. By the assumption that
$T(0)^\perp$ reduces $S$, we get that $(x,f_1)\in S$, and consequently
$(x,f_1)\in {\hat S_s}\subset S_s$ by (3.20). Thus, $(0,f_2)\in S_s$ by (3.22).
This yields that $f_2=0$ because $S_s$ is an operator. Hence,
$(x,f)=(x,f_1)\in {\hat S_s}$, and so (3.21) holds. Therefore,
the first relation in (3.19) holds.

Next, we show that
$$D\subset D({\hat A_s}).                                                \eqno(3.23)  $$
In fact, for any given $x\in D$, $x\in D(A)\cap T(0)^\perp$
by (3.15) and the assumption that $D\subset D(A)$. So there exists
$f\in X$ such that $(x,f)\in A$. Further, there exit
$f_1\in T(0)^\perp$ and $f_2\in T(0)$ such that
$(x,f)=(x,f_1)+(0,f_2)$. By the assumption that
$T(0)^\perp$ reduces $A$ we have that $(x,f_1)\in A$
and so $(x,f_1)\in {\hat A_s}$. Thus, $x\in D({\hat A_s})$,
and consequently (3.23) holds.

Finally, we show that (3.18) holds. We first show that
$$T_sx={\hat S_s}x+{\hat A_s}x,\;\;x\in D.                                \eqno(3.24) $$
It follows from (3.1), (3.15) and (3.20) that
$${\hat S_s}+{\hat A_s}\subset (S_s+A_s)\cap (T(0)^\perp)^2
\subset T\cap (T(0)^\perp)^2=T_s.                                        \eqno(3.25) $$
By (3.15), the first relation in (3.19), and (3.23) we
get that
$$D(T_s)=D=D(S_s)=D({\hat S_s})=D({\hat S_s}+{\hat A_s}).$$
Hence, by the fact that ${\hat S_s}$, ${\hat A_s}$ and $T_s$ are
(single-valued) operators, (3.25) yields (3.24),
which, together with (3.19) and
(3.23), implies that (3.18) holds.
The proof is complete.
\medskip

\noindent{\bf Remark 3.1.} If $T$ and $S$ are self-adjoint subspaces in $X^2$,
then $T_s$ and $S_s$ are self-adjoint
operators in $T(0)^{\perp}$ and $S(0)^{\perp}$, respectively, by Lemma 2.2.
So the various spectra of subspace $S$ are only affected by $ A_s|_{D(S)}$ by Lemma 2.1,
Theorem 3.3 and Theorem 4.1 of [13]
under the assumptions of Theorem 3.3.\medskip

Now, we are ready to give out relationships
between compact perturbation of a closed Hermitian subspace
and relatively compact perturbation of the operator parts
of the subspace and its perturbation term, one of which
is a characterization of compact perturbation of the subspaces
in terms of their operator parts under
 the assumptions of Theorem 3.3. \medskip

\noindent{\bf Theorem 3.4.} Let $T$, $S$ and $A$ be closed Hermitian
subspaces in $X^2$ with $D(T)=D(S)=:D\subset D(A)$, and satisfy (3.1).
Assume that $\rho(T)\cap \rho(S)\ne \emptyset$.
\begin{itemize}
\item[{\rm (i)}] If $A_s$ is an $S_s$-compact operator,
then $T$ is a compact perturbation of $S$.

\item[{\rm (ii)}] If $T$ is a compact perturbation of $S$,
$T(0)^\perp$ reduces $S$ and $A$,
and $T_s$ is a bounded operator in $D$, then $A_s$ is
an $S_s$-compact operator.

\end{itemize}
\medskip

\noindent{\bf Proof.} By the assumptions, (3.15) holds.
For convenience, denote
 $$B:=(T-\ld I)^{-1}-(S-\ld I)^{-1}.                                    \eqno(3.26)$$

 (i) Suppose that $A_s$ is an $S_s$-compact operator.
In order to show that $T$ is a compact perturbation of $S$,
it suffices to show that $B$ is a compact operator
for all $\ld \in \rho(T)\cap \rho(S)$ by Lemma 2.5.
Fix any $\ld \in \rho(T)\cap \rho(S)$ and
any bounded sequence $\{f_n\}_{n=1}^\infty\subset X$.
We want to show that $\{Bf_n\}$ has a convergent subsequence.

Set
$$x_n=(T-\ld I)^{-1}f_n,\;\;y_n=(S-\ld I)^{-1}f_n,\;\;n\ge 1.$$
Then $\{x_n\}$ and $\{y_n\}$ are bounded sequences in $D$,
$$Bf_n=x_n-y_n,\;\;n\ge 1,                                                 \eqno(3.27)$$
and
$$(x_n,f_n+\ld x_n)\in T,\;\; (y_n,f_n+\ld y_n)\in S,\;\;n\ge 1.            \eqno(3.28)$$
By noting that $y_n\in D\subset S(0)^\perp$ by (3.15),
it follows from (2.2) that there exist $g_n\in S(0)^\perp$
and $h_n\in S(0)$ such that
$$f_n=g_n+h_n,                                $$
and
$$(y_n,g_n+\ld y_n)\in S_s,\;\;(0,h_n)\in S_\infty.$$
Since $\{f_n\}$ is bounded, $\{g_n\}$ is bounded, and
consequently $\{g_n+\ld y_n\}$ is bounded.
Thus, by the assumption that $A_s$ is an $S_s$-compact operator
we have that $\{A_sy_n\}$ has a convergent subsequence.
For simplicity, denote $u_n=A_sy_n$ for $n\ge 1$.
It is evident that
$$(y_n, f_n+\ld y_n+u_n)\in S+A_s\subset T,$$
which, together with the first relation in (3.28), implies that
$(y_n-x_n, \ld(y_n-x_n)+u_n)\in T$. This yields that
$(y_n-x_n, u_n)\in T-\ld I,$
and then
$$y_n-x_n=(T-\ld I)^{-1}u_n.$$
Hence, $\{y_n-x_n\}$, namely, $\{Bf_n\}$ by (3.27), has a convergent
subsequence since $(T-\ld I)^{-1}$ is a bounded
operator on $X$. Therefore, $B$ is a compact operator in $X$,
and then $T$ is a compact perturbation of $S$.

(ii) Now, suppose that $T$ is a compact perturbation of $S$,
$T(0)^\perp$ reduces $S$ and $A$, and $T_s$ is a bounded operator.
Then $B$ is a compact operator in $X$ by Lemma 2.5

 Since all the assumptions of Theorem 3.3 hold,
$T_s$, $S_s$, $A_s$, ${\hat S_s}$,
and ${\hat A_s}$ satisfy (3.18) and (3.19),
where ${\hat A_s}$ and ${\hat S_s}$ are defined by (3.19).
Note that $\rho(T)\cap \rho(S)=\rho(T_s)\cap \rho(S_s)$ by Lemma 2.1.
So it can be easily verified by (3.18) that
for any $\ld\in \rho(T)\cap \rho(S)$,
$$(S_s-\ld I)^{-1}g=(T_s-\ld I)^{-1}g+
(T_s-\ld I)^{-1} A_s(S_s-\ld I)^{-1}g,\;\;g\in T(0)^{\perp}.     \eqno(3.29) $$
In addition, it follows from (2.2) that
$$(T-\ld I)^{-1}=(T_s-\ld I)^{-1}\oplus T_\infty^{-1},\;\;
(S-\ld I)^{-1}=(S_s-\ld I)^{-1}\oplus S_\infty^{-1}.                  \eqno(3.30) $$
For any given $f\in X$, there exist $g\in T(0)^\perp$
and $h\in T(0)$ such that $f=g+h$. By (3.30) one has
$$(T_s-\ld I)^{-1}g=(T-\ld I)^{-1}f,\;
 (S_s-\ld I)^{-1}g=(S-\ld I)^{-1}f,                                  \eqno(3.31)  $$
where the first relation in (3.19) and the assumption
that $T(0)^\perp$ reduces $S$ have been used for
the above second relation.
So, it follows from (3.26), (3.29) and (3.31) that
$$Bf=(T_s-\ld I)^{-1}g-(S_s-\ld I)^{-1}g
=-(T_s-\ld I)^{-1}A_s(S_s-\ld I)^{-1}g.
                                                                     \eqno(3.32) $$

 Fix any bounded sequence
$\{(y_n,S_sy_n)\}_{n=1}^\infty\subset S_s$.
Set $g_n=(S_s-\ld I)y_n$. Then $g_n\in T(0)^\perp$
by the first relation in (3.19), $\{g_n\}$ is a bounded sequence,
and $y_n=(S_s-\ld I)^{-1}g_n$. By (3.32) we get that
$$Bg_n=-(T_s-\ld I)^{-1} A_s y_n,\;\;n\ge 1,     $$
and then
$$A_s y_n=-(T_s-\ld I)Bg_n,\;\;n\ge 1.                               \eqno(3.33)$$
Since $B$ is compact, there exists a subsequence $\{g_{n_k}\}$
such that $\{Bg_{n_k}\}$ is convergent.
Thus,  $\{(T_s-\ld I)Bg_{n_k}\}$ is convergent by the assumption
that $T_s$ is bounded, and consequently
so is $\{A_sy_{n_k}\}$ by (3.33). Therefore, $A_s$ is
an $S_s$-compact operator. The whole proof is complete.\medskip

\noindent{\bf Remark 3.2.} By Theorem 3.4, one can see that
the relatively compact perturbation of operator parts of
closed Hermitian subspaces is stronger
than the compact perturbation, in general.\medskip

If $S$ is a self-adjoint subspace in $X^2$,
then we can give the following results:\medskip

\noindent{\bf Theorem 3.5.} Let $T$ and $A$ be closed Hermitian
subspaces in $X^2$, $S$ a self-adjoint subspace in $X^2$
with $D(T)=D(S)=:D\subset D(A)$, and they satisfy (3.1).
Assume that $\rho(T)\cap \rho(S)\ne \emptyset$.
\begin{itemize}
\item[{\rm (i)}] If $QA_s$ is an $S_s$-compact operator,
then $T$ is a compact perturbation of $S$.

\item[{\rm (ii)}] If $T$ is a compact perturbation of $S$
and $T_s$ is a bounded operator in $D$, then $QA_s$ is
an $S_s$-compact operator.

\end{itemize}
Here $Q$ is specified in Theorem 3.2.
\medskip

\noindent{\bf Proof.} The proof of assertion (ii) is similar to that of
(ii) of Theorem 3.4, where (3.18) is replaced by (3.16). So we omit
its details.

Now, we show that assertion (i) holds. Suppose that
$QA_s$ is an $S_s$-compact operator. Since all the assumptions
of Theorem 3.2 are satisfied, (3.15), (3.16) and (3.17) hold.
Fix any $\ld \in \rho(T)\cap \rho(S)$. Then $\ld \in \rho(T_s)\cap \rho(S_s)$
by Lemma 2.1. For any given $f\in X$, set
$x=(T-\ld I)^{-1}f$ and $y=(S-\ld I)^{-1}f$.
Then $(x,f+\ld x)\in T$ and $(x, f+\ld y)\in S$.
There exists $g\in T(0)^\perp$ and $h\in T(0)$
such that $f=g+h$. So we get that
$$(x,f+\ld x)=(x,g+\ld x)+(0,h),\;(y,f+\ld y)=(y,g+\ld y)+(0,h).$$
Note that $x\in D\subset T(0)^\perp$ and $S(0)=T(0)$ by (3.17).
Hence, $(x,g+\ld x)\in T_s$ and $(y,g+\ld y)\in S_s$ by (3.15),
which implies that
$$x=(T_s-\ld I)^{-1}g,\;\;y=(S_s-\ld I)^{-1}g.$$
This shows that
$$Bf=(T_s-\ld I)^{-1}g-(S_s-\ld I)^{-1}g,                       \eqno(3.34) $$
where $B$ is defined by (3.26).
On the other hand, it follows from (3.16) that
$$(S_s-\ld I)^{-1}g'=(T_s-\ld I)^{-1}g'+
(T_s-\ld I)^{-1} Q A_s(S_s-\ld I)^{-1}g',\;\;\forall\,g'\in T(0)^{\perp},      $$
which, together with (3.34), implies that
$$Bf=-(T_s-\ld I)^{-1} QA_s(S_s-\ld I)^{-1}g.                        \eqno(3.35)$$

Now, fix any bounded sequence $\{f_n\}_{n=1}^\infty\subset X$.
We shall show that $\{Bf_n\}$ has a convergent subsequence.
There exist $g_n\in T(0)^\perp$
and $h_n\in T(0)$ such that
$$f_n=g_n+h_n,   \;\;n\ge 1.                             $$
Then $\{g_n\}$ is bounded. By (3.35) one has that
$$Bf_n=-(T_s-\ld I)^{-1} QA_s(S_s-\ld I)^{-1}g_n.                        \eqno(3.36)$$
Set $y_n = (S_s-\ld I)^{-1}g_n$ for $n\ge 1$. Then
we get from (3.36) that
$$Bf_n=-(T_s-\ld I)^{-1}Q A_sy_n,\;\;n\ge 1,                      \eqno(3.37)$$
and $\{(y_n, g_n+\ld y_n)\}\subset S_s$ is bounded by the fact
that $(S_s-\ld I)^{-1}$ is a bounded operator. Thus,
 $\{QA_sy_n\}$ has a convergent subsequence by the assumption.
And consequently, $\{Bf_n\}$ has a convergent subsequence.
This means that $B$ is a compact operator. Therefore,
$T$ is a compact perturbation of $S$ by Lemma 2.5.
The proof is complete.\medskip

We shall remark that the results of Theorem 3.5 can not be directly
derived from Theorem 3.4. It is evident that the results of Theorem 3.5
are better than those of Theorem 3.4 in the case that
$S$ is a self-adjoint subspace in $X^2$.

To end this section, we give relationships
between finite rank perturbation of a closed Hermitian subspace
and finite rank perturbation of the operator parts
of the subspace and its perturbation term, one of which
is a characterization of finite rank perturbation of the subspace
in terms of their operator parts.\medskip

\noindent{\bf Theorem 3.6.} Let $T$, $S$ and $A$ satisfy the assumptions
of Theorem 3.4.
\begin{itemize}
\item[{\rm (i)}] If $A_s$ is a finite rank operator in $D$,
then $T$ is a finite rank perturbation of $S$.

\item[{\rm (ii)}] If $T$ is a finite rank perturbation of $S$ and
$T(0)^\perp$ reduces $S$ and $A$,
then $A_s$ is a finite rank operator in $D$.
\end{itemize}
\medskip

\noindent{\bf Proof.} With a similar argument to that used in the proof
of Theorem 3.4, one can easily show Theorem 3.6 by Lemmas 2.1 and 2.5,
Theorem 3.3 and (3.32). So we omit its details.\medskip

Note that the assumption that $T_s$ is a bounded operator in $D$
in (ii) of Theorem 3.4 has been removed in (ii) of Theorem 3.6
because that if $B$ is finite rank, then $(T_s-\ld I)B$ is finite
rank for every linear operator $T_s$.\medskip

\noindent{\bf Theorem 3.7.} Let $T$, $S$ and $A$ satisfy the assumptions
of Theorem 3.5. Then $QA_s$ is a finite rank operator in $D$
if and only if $T$ is a finite rank perturbation of $S$,
where $Q$ is specified in Theorem 3.2.
\medskip

\noindent{\bf Proof.} With a similar argument to that used in the proof
of Theorem 3.5, one can easily show Theorem 3.7 by Lemmas 2.1 and 2.5,
Theorem 3.2 and (3.35). So we omit its details.\medskip

\noindent{\bf Remark 3.3.} The finite rank property of $QA_s$
is equivalent to the finite rank perturbation of the self-adjoint subspace $S$
under the assumptions of Theorem 3.7. It is evident that
if $A_s$ is finite rank in $D=D(S)$, then so is $QA_s$.
Its converse may not hold in general. But, in the case that
$S(0)\ominus A(0)$ is finite-dimensional, then the converse is true because
$A_s=QA_s+(I-Q)A_s$, while $I-Q:\,A(0)^\perp\to A(0)^\perp\ominus S(0)^\perp$
is finite rank by the fact that $A(0)^\perp\ominus S(0)^\perp=S(0)\ominus A(0)$. \medskip

\noindent{\bf Remark 3.4.} A trace class perturbation is a special
compact perturbation and a finite rank perturbation is a simple
case of trace class perturbation. We shall study the trace class perturbation of
closed subspaces in detail in our forthcoming paper.
\bigskip

\noindent{\bf 4. Self-adjoint subspaces under compact perturbations}\medskip

In this section, we show that a self-adjoint subspace is still self-adjoint
under relatively bounded and relatively compact perturbations as well as compact perturbation.\medskip

The following result is a generalization of the well-known
Kato-Rellich theorem for self-adjoint operators (cf.,
[9, Theorem 10.12] or [11, Theorem 8.5]
or [16, Theorem 5.28]) to self-adjoint subspaces.\medskip

\noindent{\bf Theorem 4.1.} Let $S$ be a self-adjoint
subspace in $X^2$ and $A$ a closed Hermitian subspace in $X^2$ with
$D(S)\subset D(A)$. If $A_s$ is $S_s$-bounded with $S_s$-bound
less than $1$, then
$S+A$ is a self-adjoint subspace in $X^2$.\medskip

\noindent{\bf Proof.} We shall show that $S+A$ is self-adjoint
in $X^2$ by Lemma 2.3.

It follows from Lemmas 2.1 and 2.2 that
 $S_s$ is a self-adjoint operator
in $S(0)^\perp$, and $A_s$ is a closed Hermitian operator in
$A(0)^\perp$ with $D(S_s)=D(S)\subset D(A)=D(A_s)$.
Then $D(S)$ is dense in $S(0)^\perp$, and consequently
$A_s$ is densely defined in $S(0)^\perp$.
Since $A_s$ is $S_s$-bounded with $S_s$-bound less than $1$,
$S_s+A_s$ is a self-adjoint operator in $S(0)^\perp$
by the Kato-Rellich theorem [9, Theorem 10.12]. Hence,
by Lemma 2.3 one has
$$R(S_s+A_s\pm iI)=S(0)^\perp.                                     \eqno(4.1)$$
Fix any $f\in X$. There exist $f_1\in S(0)^\perp$ and $f_2\in S(0)$
such that $f=f_1+f_2$.
Further, by (4.1), there exists $x\in D(S_s)$ such that
$$f_1=S_sx+A_sx\pm ix.                                                      \eqno(4.2) $$
Note that $(x,S_sx)\in S_s\subset S$, $(x,A_sx)\in A_s\subset A$,
and $(0,f_2)\in S_\infty\subset S$, which implies that
$(x, S_sx+A_sx+f_2)\in S+A$. This, together with (4.2), yields that
$(x,f)=(x, S_sx+A_sx+f_2\pm ix)\in S+A\pm iI$. This means
that $f\in R(S+A\pm iI)$. Hence, $R(S+A\pm iI)=X$, and
consequently $S+A$ is a self-adjoint subspace in $X^2$ by Lemma 2.3.
 This completes the proof.

\medskip

If the condition that $S+A$ is closed is added to the assumptions
of Theorem 4.1, then the condition that
$A_s$ is $S_s$-bounded with $S_s$-bound
less than $1$ can be weakened as follows.
We shall point out again that $S+A$ may not be closed
if $S$ and $A$ are closed.\medskip

\noindent{\bf Theorem 4.2.} Let $S$ be a self-adjoint
subspace in $X^2$, and $A$ and $S+A$ closed Hermitian subspaces in $X^2$ with
$D(S)\subset D(A)$. If $QA_s$ is $S_s$-bounded with $S_s$-bound
less than $1$, then $S+A$ is a self-adjoint subspace in $X^2$,
where $Q$ is specified in Theorem 3.2.\medskip

\noindent{\bf Proof.} Let $T$ be defined by (3.1).
By Theorem 3.2, $T_s$, $S_s$ and $A_s$ satisfy (3.16) and (3.17) holds. Further, $S_s$ is a self-adjoint operator
in $S(0)^\perp$ by Lemma 2.2, and  $QA_s$ is an Hermitian operator
in $D(S)\subset D(A)$ by Theorem 3.1.
So  $QA_s$ is densely defined and Hermitian in $S(0)^\perp$,
where the fact that $D(S)$ is dense in $S(0)^\perp$ has been used.
Since $QA_s$ is $S_s$-bounded with $S_s$-bound
less than $1$ by the assumption, it follows from (3.16) that
$T_s$ is a self-adjoint operator in $S(0)^\perp$ by
the Kato-Rellich theorem [9, Theorem 10.12].
Note that $T(0)^\perp=S(0)^\perp$ by (3.17).
Therefore, $T$, namely, $S+A$, is a self-adjoint subspace in $X^2$
by Proposition 2.1. This completes the proof.\medskip

By Theorem 9.7 of [16], if an operator $U$ is relatively compact
to another operator $V$, then $U$ is
$V$-bounded with $V$-bound zero. So the following
result can be derived from Theorem 4.1.\medskip

\noindent{\bf Theorem 4.3.}  Let $S$ be a self-adjoint
subspace in $X^2$ and $A$ a closed Hermitian subspace in $X^2$ with
$D(S)\subset D(A)$. If $A_s$ is $S_s$-compact, then
$S+A$ is a self-adjoint subspace in $X^2$.\medskip

Since a finite rank or more general trace class
operator is compact, the following result
can be directly derived from Theorem 4.3.\medskip

\noindent{\bf Corollary 4.1.}  Let $S$ be a self-adjoint
subspace in $X^2$ and $A$ a closed Hermitian subspace in $X^2$ with
$D(S)\subset D(A)$. If $A_s$ is finite rank or belongs
to the trace class in $D(S)$, then
$S+A$ is a self-adjoint subspace in $X^2$.\medskip

The following result is a direct consequence of Theorem 4.2.\medskip

\noindent{\bf Theorem 4.4.} Let $S$ be a self-adjoint
subspace in $X^2$, and $A$ and $S+A$ closed Hermitian subspaces in $X^2$ with
$D(S)\subset D(A)$. If $QA_s$ is $S_s$-compact,
then $S+A$ is a self-adjoint subspace in $X^2$,
where $Q$ is specified in Theorem 3.2.\medskip

The following result can be directly derived from
Theorem 4.4.\medskip

\noindent{\bf Corollary 4.2.} Let $S$ be a self-adjoint
subspace in $X^2$, and $A$ and $S+A$ closed Hermitian subspaces in $X^2$ with
$D(S)\subset D(A)$. If $QA_s$ is  finite rank or belongs
to the trace class in $D(S)$,
then $S+A$ is a self-adjoint subspace in $X^2$,
where $Q$ is specified in Theorem 3.2.\medskip

The following result is a direct consequence of Theorems 3.5 and 4.4.\medskip

\noindent{\bf Corollary 4.3.} Let $S$ be a self-adjoint subspace in $X^2$,
and $A$ and $S+A$ closed Hermitian subspaces in $X^2$ with
$D(S)\subset D(A)$. Assume that $\rho(T)\cap \rho(S)\ne \emptyset$.
If $S+A$ is a compact perturbation of $S$
and $(S+A)_s$ is a bounded operator, then $S+A$ is a self-adjoint subspace in $X^2$.

\medskip

If we strengthen the assumption about
the resolvent sets of $S$ and $S+A$ in Corollary 4.3,
the assumption that $(S+A)_s$ is a bounded operator can be removed as follows.
\medskip

\noindent{\bf Theorem 4.5.}  Let $S$ be a self-adjoint
subspace in $X^2$ and $T$ a closed
Hermitian subspace in $X^2$.
Assume that $\rho(S)\cap \rho(T)\cap {\mathbf R} \ne \emptyset$.
If  $T$ is a compact perturbation of $S$,
then $T$ is a self-adjoint subspace $X^2$.\medskip

\noindent{\bf Proof.} Fix any $\ld \in \rho(S)\cap \rho(T)
\cap {\mathbf R}$. By the definition of self-adjoint subspace,
it can be easily verified that
a subspace $C$  is self-adjoint (resp., closed Hermitian) in $X^2$
if and only if $C^{-1}$ is self-adjoint (resp., closed Hermitian)
in $X^2$. Hence, $(S-\ld I)^{-1}$ is a bounded self-adjoint
operator on $X$ and $(T-\ld I)^{-1}$ is a bounded and closed
Hermitian operator on $X$. Let $B$ be defined by (3.26).
Then $B$ is a compact Hermitian operator defined on $X$
and
 $$(T-\ld I)^{-1}=(S-\ld I)^{-1}+B,$$
which yields that $(T-\ld I)^{-1}$ is a self-adjoint operator on $X$
by Lemma 2.4. Thus, $T-\ld I$, and then $T$ by (2.1),
is a self-adjoint subspace in $X$.
This completes the proof.\medskip

The following result directly follows from Theorem 4.5.\medskip

\noindent{\bf Corollary 4.4.}  Let $S$ be a self-adjoint
subspace in $X^2$ and $T$ a closed
Hermitian subspace in $X^2$.
Assume that $\rho(S)\cap \rho(T)\cap {\mathbf R} \ne \emptyset$.
If $T$ is a finite rank perturbation of $S$
or $(T-\ld I)^{-1}-(S-\ld I)^{-1}$ belongs to the trace class in $X$
for some $\ld \in \rho(S)\cap \rho(T)\cap {\mathbf R}$,
then $T$ is a self-adjoint subspace in $X^2$.\medskip

\noindent{\bf Remark 4.1.} The finite rank perturbation
is of special interest in the existing literature.
So we specially point out this condition in the above results.

\bigskip

\noindent{\bf 5. Stability of essential spectra
of self-adjoint subspaces under perturbations}\medskip

In this section, we study stability of essential spectra
of self-adjoint subspaces under compact and relatively compact perturbations.\medskip

\noindent{\bf Theorem 5.1.}  Let $T$ and $S$  be self-adjoint subspaces
in $X^2$. If $T$ is a compact perturbation of $S$,
then
$$\si_e(T)=\si_e(S),                                           \eqno(5.1) $$
and consequently $S$ has a pure discrete spectrum if and only
if so does $T$.\medskip

\noindent{\bf Proof.} By Definition 2.4,
$T$ is a compact perturbation of $S$ if and only if
$S$ is a compact perturbation of $T$.
So it suffices to show that
$$\si_e(S)\subset \si_e(T)                                      \eqno(5.2)$$
because its inverse inclusion can be obtained by interchanging
$S$ and $T$.

Fix any $\ld \in \si_e(S)$. Then, by Theorem 3.7 of [13],
there  exists a  sequence
$\{(x_n,f_n)\} \subset S$ satisfying that
$$x_n \stackrel{w}{\to} 0,\;\; \liminf_{n\to \infty}\|x_n\|> 0,\;\;
 f_n-\ld x_n \to 0\;{\rm as}\;n\to \infty.                                 \eqno(5.3)$$
Since $T$ and $S$  are  self-adjoint subspaces in $X^2$,
${\mathbf C}\setminus {\mathbf R}\subset \rho(T)\cap \rho(S)$
by Theorem 2.5 of [13], and then
$\rho(T)\cap \rho(S)\ne \emptyset$.
Take any $\mu\in \rho(T)\cap \rho(S)$. Then
$B$ is a compact operator defined on $X$ by Lemma 2.5,
where $B$ is defined by (3.26) with $\ld$ replaced by $\mu$.
Further, noting that $(x_n,f_n-\mu x_n)\in S-\mu I$, we have
that
$$x_n=(S-\mu I)^{-1}(f_n-\mu x_n).                                       \eqno(5.4) $$
Set
$$y_n=(T-\mu I)^{-1}(f_n-\mu x_n).                                       \eqno(5.5)$$
Then $(y_n, f_n-\mu x_n)\in T-\mu I$, which implies that
$(y_n, f_n+\mu (y_n- x_n))\in T$. Now, we want to show that
$$y_n \stackrel{w}{\to} 0,\;\; \liminf_{n\to \infty}\|y_n\|> 0,\;\;
 f_n+\mu (y_n- x_n)-\ld y_n \to 0\;{\rm as}\;n\to \infty.                  \eqno(5.6)$$
If it is true, then $\ld\in \si_e(T)$ again by Theorem 3.7 of [13],
and consequently, (5.2) holds. It follows from (5.4) and (5.5) that
$$y_n-x_n=B(f_n-\mu x_n).                                           \eqno(5.7) $$
In addition, by (5.3) we get that
$$f_n-\mu x_n=f_n-\ld x_n+(\ld-\mu) x_n\stackrel{w}{\to} 0,$$
which, together with (5.7), the compactness of $B$
and Theorem 6.3 of [16], yields that
$$y_n-x_n\to 0\;{\rm as}\;n\to \infty.                                    \eqno(5.8) $$
This implies that the first relation in
(5.6) holds, and
$\liminf_{n\to \infty}\|y_n\|=\liminf_{n\to \infty}\|x_n\|>0$.
Moreover, we have
$$f_n+\mu (y_n- x_n)-\ld y_n=f_n-\ld x_n+(\mu-\ld) (y_n- x_n),$$
which, together with the third relation in (5.3) and (5.8), implies that
the third relation in (5.6) holds. Hence, (5.6) has been shown.
Therefore, (5.1) holds.

The final assertion can be directly derived from (5.1)
and (4) of Definition 2.1. The proof is complete.\medskip

By Theorems 4.5 and 5.1 one can easily get the following result:\medskip

\noindent{\bf Corollary 5.1.}  Let $S$ be a self-adjoint
subspace in $X^2$ and $T$ a closed
Hermitian subspace in $X^2$.
Assume that $\rho(S)\cap \rho(T)\cap {\mathbf R} \ne \emptyset$.
If $T$ is a compact perturbation of $S$,
then the results of Theorem 5.1 hold.\medskip

The following result is a direct consequence of Corollary 5.1.\medskip

\noindent{\bf Corollary 5.2.}  Let $S$ be a self-adjoint
subspace in $X^2$ and $T$ a closed
Hermitian subspace in $X^2$.
Assume that $\rho(S)\cap \rho(T)\cap {\mathbf R} \ne \emptyset$.
If $T$ is a finite rank perturbation of $S$
or $(T-\ld I)^{-1}-(S-\ld I)^{-1}$ belongs to the trace class in $X$
for some $\ld \in \rho(S)\cap \rho(T)\cap {\mathbf R}$,
then the results of Theorem 5.1 hold.\medskip

The following result can be directly derived from Theorem 5.1. \medskip

\noindent{\bf Corollary 5.3.} Let $T$ and $S$  be self-adjoint subspaces
in $X^2$. If $T$ is a finite rank perturbation of $S$
or $(T-\ld I)^{-1}-(S-\ld I)^{-1}$  belongs to the trace class in $X$
for some $\ld \in \rho(S)\cap \rho(T)$,
then the results of Theorem 5.1 hold.\medskip

\noindent{\bf Theorem 5.2.}  Let $S$  be a self-adjoint subspace
in $X^2$ and $A$ a closed Hermitian subspace in $X^2$ with
$D(S)\subset D(A)$. If $A_s$ is an $S_s$-compact
operator, then
$$\si_e(S)= \si_e(S+A),                                              \eqno(5.9)$$
and consequently $S$ has a pure discrete spectrum if and only
if so does $S+A$.\medskip

\noindent{\bf Proof.} Let $T$ be defined by (3.1).
 By Theorem 4.3, $T$ is a self-adjoint subspace
in $X^2$. Thus, ${\mathbf C}\setminus {\mathbf R}
\subset \rho(T)\cap \rho(S) $ by Theorem 2.5 of [13],
and consequently $\rho(T)\cap \rho(S)\ne \emptyset$.
Further, it follows from Theorem 3.4 that
$T$ is a compact perturbation of $S$.
Therefore, (5.9) holds by Theorem 5.1.
This completes the proof.\medskip

The following result directly follows from Theorem 5.2.\medskip

\noindent{\bf Corollary 5.4.}  Let $S$  be a self-adjoint subspace
in $X^2$ and $A$ a closed Hermitian subspace in $X^2$ with
$D(S)\subset D(A)$. If $A_s$ is finite rank or belongs to
the trace class in $D(S)$, then the results of Theorem 5.2 hold.\medskip

\noindent{\bf Theorem 5.3.} Let $S$ be a self-adjoint
subspace in $X^2$, and $A$ and $S+A$ closed Hermitian subspaces in $X^2$ with
$D(S)\subset D(A)$. If $QA_s$ is $S_s$-compact,
then the results of Theorem 5.2 hold,
where $Q$ is specified in Theorem 3.2.\medskip

\noindent{\bf Proof.} By Theorem 4.4, $T=S+A$ is a self-adjoint subspace
in $X^2$, and by Theorem 3.2, $T_s$, $S_s$ and $A_s$ satisfy (3.16) and
(3.17) holds. So $T_s$ and $S_s$ are self-adjoint operators
in $S(0)^\perp=T(0)^\perp$ by Lemma 2.2 and (3.17). Further, $QA_s$
is an Hermitian operator in $D(S)$ by Theorem 3.1.
By noting that $D(S)$ is dense in $S(0)^\perp$,
$QA_s|_{D(S)}$ is a densely defined and Hermitian operator in $S(0)^\perp$.
Therefore, $\si_e(T_s)=\si_e(S_s)$ by Lemma 2.4, (3.16) and the assumption
that $QA_s$ is $S_s$-compact. This yields the results of Theorem 5.3
by Lemma 2.1, and then the proof is complete.\medskip

By Theorem 5.3 one can easily get the following result:\medskip

 \noindent{\bf Corollary 5.5.} Let $S$ be a self-adjoint
subspace in $X^2$, and $A$ and $S+A$ closed Hermitian subspaces in $X^2$ with
$D(S)\subset D(A)$. If $QA_s$ is finite rank or belongs to
the trace class in $D(S)$, then the results of Theorem 5.2 hold,
where $Q$ is specified in Theorem 3.2.\medskip

Finally, we shall give several remarks on the assumptions and results
of Theorems 5.1-5.3 and Corollaries 5.1 - 5.5.\medskip

\noindent{\bf Remark 5.1.} Theorems 5.1 and 5.2 extend Theorems 8.12
and 8.15 of [11] for self-adjoint operators to self-adjoint subspaces,
respectively.\medskip

\noindent{\bf Remark 5.2.} In the case
that it is known that a subspace and its perturbation are both
self-adjoint, and their resolvents can be explicitly expressed,
then Theorem 5.1 and Corollary 5.3 are applicable.
In the case that it is known that the resolvent sets
of the unperturbed self-adjoint subspace and its perturbation
both contain a real value, then Corollaries 5.1 and 5.2
are applicable. Instead, if it is more easy to get the related
information of the operator part of the perturbed term,
then Theorems 5.2 and 5.3 and Corollaries 5.4 and 5.5 are more applicable.
\medskip

\noindent{\bf Remark 5.3.} As we have mentioned in the first section,
it is very important for us to study spectral properties of
multi-valued or non-densely defined Hermitian
operators because a minimal operator, generated by a symmetric linear
difference or differential expression that does not
satisfy the corresponding definiteness condition,
may be multi-valued and non-densely defined
and so may be its self-adjoint extensions
(cf., [8], [10] and [14]).
The results given in this section are available in this case. \medskip

\noindent{\bf Remark 5.4.} By Theorems 5.1 and 5.2, the essential spectrum
$\si_e(S)$ of a self-adjoint subspace $S$ is invariant if
the perturbation is compact or the operator part
$A_s$ of the perturbed term $A$ is $S_s$-compact.
if the perturbation or
$A_s$ is finite rank or more generally  belongs to the trace class,
then $\si_e(S)$ is invariant by Corollaries 5.2 - 5.5.
We shall further study invariance of the absolutely continuous
spectrum of $S$ under this perturbation in our forthcoming paper.
In addition, we shall apply these results
to study dependence of absolutely continuous spectra
on regular endpoints and boundary conditions
and invariance of essential and absolutely continuous spectra
under perturbation for
singular linear Hamiltonian systems in our other
forthcoming papers, including that the systems are in
the limit point and middle limit cases
at singular endpoints.\bigskip

\noindent{\bf Acknowledgements}\medskip

This research was partially supported by the NNSF of Shandong Province (Grant ZR2011
\newline AM002).

This paper was completed when the author visited California Institute of Technology, U.S.A.,
during November 2012 - May 2013. She would like to thank the China Scholarship Council, P. R. China, for financial support, and thank Professor Barry Simon of the Department of Mathematics, California Institute of Technology, for hosting her visit with excellent research support.

\bigskip

\bigskip \noindent{\bf \large References}
\def\hang{\hangindent\parindent}
\def\textindent#1{\indent\llap{#1\enspace}\ignorespaces}
\def\re{\par\hang\textindent}
\noindent \vskip 3mm

\re{[1]} R. Arens, Operational calculus of linear relations,
Pac. J. Math. 11(1961) 9--23.

\re{[2]} T. Ya. Azizov, J. Behrndt, P. Jonas, C. Trunk,
Compact and finite rank perturbations of linear
relations in Hilbert spaces, Integr. Equ. Oper. Theory 63 (2009) 151-163.

\re{[3]} T. Ya. Azizov, J. Behrndt, P. Jonas, C. Trunk, Spectral points of definite type
and type $\pi$ for linear operators and relations in Krein spaces,
J. Lond. Math. Soc. 83 (2011) 768--788.

\re{[4]} A. Dijksma, H. S. V. de Snoo,
Eigenfunction extensions associated with pairs of ordinary differential expressions,
J. Differ. Equations 60(1985) 21--56.

\re{[5]} S. Hassi, H. S. V. de Snoo, One-dimensional graph
perturbations of self-adjoint relations, Ann. Aca. Sci. Fenn.
Math. 20(1997) 123--164.

\re{[6]} S. Hassi, H. de Snoo, F. H. Szafraniec,
Componentwise and cartesian decompositions of linear relations,
Dissertationes Math. 465, 2009 (59 pages).

\re{[7]} T. Kato, Perturbation Theory for Linear Operators, 2nd ed., Springer-Verlag, Berlin \newline/Heidelberg/New York/Tokyo,  1984.

\re{[8]} M. Lesch, M. Malamud,
On the deficiency indices and self-adjointness of
symmetric Hamiltonian systems, J. Differ. Equations 18(2003) 556--615.

\re{[9]} M. Reed, B. Simon, Methods of Modern Mathematical Physics
II: Fourier Analysis, Self-adjointness, Academic Press, 1972.

\re{[10]} G. Ren, Y. Shi,
Defect indices and definiteness conditions
for discrete linear Hamiltonian systems, Appl. Math. Comput.
218(2011) 3414--3429.

\re{[11]} K. Schm\"udgen, Unbounded self-adjoint operators on Hilbert space,
Springer Dordrecht Heidelberg/ New York/ London, 2012.


\re{[12]} B. Simon, Trace Ideals and Their Applications, Amer. Math. Soc.,
2nd ed., Providence, Rhode Island, 2005.


\re{[13]} Y. Shi, C. Shao, G. Ren, Spectral properties of self-adjoint subspaces,
Linear Algebra Appl. 438(2013) 191--218.

\re{[14]} Y. Shi, H. Sun,
Self-adjoint extensions for second-order
symmetric linear difference equations, Linear Algebra Appl.
434(2011) 903--930.

\re{[15]} G. Teschl, Jacobi Operators and Completely Integrable Nonlinear
Lattices, Amer. Math. Soc., Providence, Rhode Island, 2000.

\re{[16]} J. Weidmann, Linear Operators in Hilbert Spaces,
Graduate Texts in Math., vol.68, Springer-Verlag, New
York/Berlin/Heidelberg/Tokyo, 1980.

\end{document}